\documentclass[12pt]{article}

\usepackage{latexsym,amsfonts,amsthm,amsmath,amscd,amssymb}
\usepackage[dvips]{graphicx}

\newcommand{\lz}{\mathrm{Lz}_{3+1}^{10}}
\newcommand{\fn}{\mathrm{Fn}_{3+1}^{10}}

\begin{document}

\title{Fantappi\`e's final relativity and Lie algebra deformations}

\author{N. Ciccoli\\
Dipartimento di Matematica e Informatica\\
Universit\`a di Perugia}

\maketitle

\section{Introduction}
In a brief paper published by Fantappi\`e in 1954 (\cite{Fa54}) the rigidity of the real semisimple Lie algebra $\mathfrak{so}(4,1)$ was first proved. The purpose of this note is to provide some historical context of this work and discuss  why no further developments of this result were pursued by Italian mathematicians in those years.

\section{Biographical notes on Fantappi\`e}

Luigi Fantappi\`e was born in Viterbo in 1901. He studied mathematics in Pisa, at Scuola Normale Superiore, where he graduated under the direction of L. Bianchi in 1922. Under the guidance of Bianchi his first mathematical research was 
on algebraic topics, namely directed to the study of Riemann zeta function (\cite{Fa73}, vol. I). Soon after graduation he moved to Rome, where he was in direct contact with two of the leading mathematicians of the time: Francesco Severi and Vito Volterra. Severi was, at that time, the rising star of Italian mathematics. He had just moved from Padua to Rome and despite the fact that he was a geometer he was appointed as Professor of Analysis. His reputation was undisputed and after one year he was elected dean of the Universit\`a \lq\lq La Sapienza". Fantappi\`e was offered a position as assistant of Severi, a role that would imply taking notes of his lectures, and implicitely following his advices. At the same time, however, Fantappi\`e had a chance to contact Vito Volterra, maybe the most internationally renowned mathematician of those years, though already in the second half of his career. Among his many and widespread scientific interests Volterra pioneered functional analysis in Italy. It was the outcome of such different interactions that stimulated  Fantappi\`e in studying a new class of functionals, generalizing previous works of Volterra form the real to the complex case. 
Beginning in 1924, Fantappi\'e defined \emph{analytic functionals} and relatively rapidly his contributions to this field were recognized internationally  (see \cite{Lu82,Str} for a historical overview of this theory). As a result of his efforts he was promoted professor and appointed to Cagliari, Palermo, Florence and Bologna.
During the ensuing years his relationship with each of his mentors diverged.
 Volterra considered Fantappi\`e \emph{his prize pupil}, as Weil recalls (\cite[p. 48]{Wei}). However they had conflicting  political opinions. While Fantappi\`e was a fervent Fascist from his arrival in Rome,  Volterra was one of the leading figures of opposition against the regime. In 1932 Volterra was dismissed from the University for his refusal to take the oath of loyalty to the Fascist government. This was not enough to close scientific contacts between the two, which were however definitely interrupted after the anti-Jewish legislation in 1938;  Fantappi\`e welcomed in front of the old, Jew, professor the new law. Conversely his relation with Severi grew stronger. Severi became the leading mathematician in Italy and the one to have direct connections with the Government. 
Fantappi\'e was an ardent supporter of Severi and they developed a close friendship strengthened by their sharing similar political, religious and scientific viewpoints.  
In 1933 Fantappi\'e accepted a position as foundation professor of mathematics at a new University in Sa\^o Paulo, Brazil.
 His mission was both to establish the new mathematics Department and to strengthen relations with the large expatriate Italian community, a mission not devoid of political issues. He remained in Brazil from 1933 to 1939\footnote{see \cite{Zor} for a detailed description of Fantappi\`e's influences on Brazilian mathematics.}), founding the local mathematical newspaper, guiding the more promising students in researches centered around analytic functionals, and in promoting the ideas of the Fascist regime in meetings and conferences throughout the country.  In 1939 he left Brazil and went back to Rome where he was offered the Chair of Higher Analysis and  was appointed Vice President of the newly founded National Institute of Higher Mathematics (INDAM), directed by his former mentor Severi. During WWII he often travelled to Spain, where he developed a close connection with the University of Catalunya, and Portugal and began enlarging his research interests. Starting with 1942 he proposed a mainly philosophical theory centered on the concept of \emph{sintropy} which, in his intentions, should connect physics and biology and explain \lq\lq finalistic phenomena". In later years he also proposed a \lq\lq mathematical" proof of the existence of God and devoted some attention (in conferences and in published papers) to parapsychological phenomena (see \cite{FaCo}). Apart from this unconventional interests, from the late 1940's his attentions were drawn to the role of topological groups in physical and mathematical theories: he lectured on topological groups at Indam and published papers centered on the role of symmetries in relativity theory. This  was to be his last mathematical interest and in the summer of 1956 he suffered a fatal heart attack. Some additional details on his life can be found in the numerous obituaries (a quite complete one is \cite{Fi}), and in the introduction to the first volume of his collected works \cite{Fa73}.

\section{Fantappi\`e's final relativity}

In looking at Fantappi\`e's list of scientific publications (as, for example, recorded in \cite{Fa73}) one sees a distinct gap between 1943 and 1948. This gap may not be immediately apparent since it does not consist of a complete lack of new publications. All his works, in those years, were not mathematical but focuse on the development of \emph{sintropy} with scientific considerations showing up just in the background.

The theory of analytic functionals which had established his reputation was largely put aside. The problem of extending analytic functionals to more variables could not be easily solved by the techniques at hand in his school (as nicely explained by Struppa in \cite{Str}). However, he did acknowledge it in lecture notes for a course in Barcelona in 1948.

There were also external reasons, both social and personal, for this interruption. Rome, the city where he was living, during 1943 experienced the worst years of the WWII period. Despite propaganda to the contrary, it was clear by that time that the Axis powers had lost and in July that year the northward Allied advance had begun. Also in July the Allied bombing of Rome commenced with approximately 3000 casualties in a single day. It was going to  be but the first of more than 50 bombings on the capital, in subsequent months. Fantappi\`e was in those days in his villa just outside Viterbo, his birthplace, only 90 kilometers from Rome; it is in fact reported that he missed one of the Indam's meetings, in subsequent days, due to difficulties in reaching Rome. Viterbo was also bombed towards the end of 1943.  At the end of July the Fascist government was overthrown and Mussolini arrested. This marked the start of a long and difficult period of social unrest. By the beginning of September the new Italian government switched allegiance to the Allies and German military forces occupied the city. In October 1943, Nazis  raided the Jewish ghetto in Rome, deporting more than 1.000 in concentration camps; all Italian young males unwilling to take arms against the Allies had to hide themselves to escape immediate execution. It is no surprise that under such exceptional circumstances the University was essentially closed. Even after Rome's liberation on the fifth of  June 1944 the situation did not get any better. Severi,  a prominent symbol of deposed Fascist regime and Fantappi\'e's mentor and personal friend, was suspended from all his University positions (\cite{GoBa}). Despite F's outspoken Fascist support there is little information he was similarly sanctioned. However, he was burdened by the additional role of Indam's vicarious president until 1948 (though, according to \cite{Ro}, Severi got back to his duties, though not officially, already in 1946).

On a personal level his life went through big changes as well. The death of his beloved mother, with whom he was living, was followed quite closely by his marriage with Maria Quadrani, a family friend,  in 1946.

Being as it is, whether social, personal or scientific reasons prevailed in such turbulent years, no new papers on pure math were published by Fantappi\`e, in this five years window. A new interest, however, arose: the role of symmetries in physical systems, with a special emphasis on relativity theory. In 1949 a note  (\cite{Fa49}) appeared on \emph{Costruzione effettiva dei prodotti funzionali relativisticamente invarianti}\footnote{An effective construction of relativistically invariant functional products; {\tt MR 0044942}.}. This was the prelude to an almost complete change of direction of his scientific interest. Though he kept on assisting former students and young colleagues, in Rome as well as in Barcellona and Sa\^o Paulo, in their efforts to further develop his theory of analytic functionals, his own research focused on the role of  symmetries in physics.  His last student, G. Arcidiacono, a physicist, was the only one to follow him along this path. This change of direction was reflected also in his Indam lectures (INDAM lectures at that time were analogous to modern day PhD courses). In the troublesome academic year of 1945/1946, the first after the end of WWII, he abandoned analytic functionals, the traditional argument of his course, and lectured about  \emph{Mathematical tools of modern physics and Universe Unitary Theory}\footnote{Strumenti matematici della fisica moderna e Teoria Unitaria dell'Universo}.  In 1950/1951 he switched again and, until his death in 1956, his lecture courses were devoted to \emph{Topological Groups and their application to Physics}\footnote{Gruppi topologici e loro applicazioni fisiche} (\cite{Ro,Suc}).

Fantappi\`e's researched on symmetries in mathematical physics had to face  a rather peculiar state of the art. 

The history of Relativity Theory in Italy has been much treated in the literature (e.g. \cite{Be,Pas}) and we will only recall some key points. After a very intense debate in the 1920's between a relativistic and an antirelativistic side, by the 1930's the situation was settled in favour of the first. This was not completely independent of the fact that the work of Levi--Civita stood  in favour of Relativity Theory, with his great scientific reputation. However Relativity was  studied more in mathematics than physics departments and after Levi-Civita was discharged from the University as a consequence of Racial Laws, it started fading also from the attention of mathematicians. Therefore when Fantappi\`e understood his researches in the field this neglect of the topic was quite evident.

Lie theory was a somewhat analogous situation. It is not too excessive to say that soon after Lie's work Italy became one of the leading nations in pursuing this line of research. Following the work of Cremona and Corrado Segre, that immediately saw how fruitful such theory could be, a number of people contributed relevant results. The first to be mentioned is certainly L. Bianchi, whose treatise \cite{Bia} rapidly became a standard in years to come. He was followed and surrounded by work of Enriques, Fano, Fubini, the already mentioned Levi-Civita, Medolaghi, Segre and Vivanti, all  between the end of the 19th and the beginning of the 20th century. A special mention should be given to the work of the young  Eugenio Elia Levi, at the outbreak of the 20th century, that gave him a rapidly growing reputation  abruptly interrupted by his death in action in WWI. All this attention faded rapidly. When few years later U. Amaldi published his works concerning  infinite transformations groups he got almost no reaction from the mathematical community ( ``The work of U. Amaldi [...] left no trace. This should not surprise if we consider that also the fundamental work of Cartan on finite transformations groups was basically ignored until the 60ies.`` \cite[p. 64]{Rog} \footnote{Il lavoro di U. Amaldi [...] pass\`o praticamente inosservato. La cosa non deve stupire se si considera che anche gli importantissimi lavori di Cartan sui gruppi continui finiti non furono praticamente letti fino agli anni '60}). It was the beginning of a long eclipse, during which Lie theory received a very limited attention in the Italian School of Mathematics. It is significant that in the standard comprehensive book on Italian Mathematics between the two World Wars \cite{Pas} the name of Lie appears only twice, and always in connection to applications, either in Differential Geometry, following the work of Bianchi, or in Integrability of Mechanical Systems. 

Thus, when Fantappi\`e undertook his new research, there was scant attention given in Italy to either of these areas. Conversely, there was much activity internationally which had generated a wealth of new results.  Fantappi\`e himself was somewhat aware of this gap and was referring to his students about his limited mathematical tools in comparison to those used by some foreign colleagues (I owe this remark to a personal communication with Prof. F. Succi, one of Fantappi\`e last students).

Fantappi\`e was guided by the idea of realizing what he himself calls ``an Erlangen program for physics'', a classification of possible physical theories through their group of symmetries. With these words he ended \cite[p. 290]{Fa52}:
\begin{quote}
 From the previous, it then follows that, just as for the classification of geometries given by Klein in his Erlangen's program, also for possible \lq\lq physics", or better \lq\lq physical universes" as above specified, the most natural possible classification is the one in terms of the fundamental groups used to define \lq\lq equality". Let us remark that here the term \lq\lq fundamental group" has to be understood as maximal group of symmetries of the physical system.\footnote{"Da quanto precede, infine, segue allora che, come gi\`a per la classificazione delle geometrie data da Klein nel suo programma di Erlangen, anche per le \lq\lq fisiche" possibili, o meglio per gli \lq\lq universi fisici" possibili, sopra specificati, la classificazione pi\'u naturale che si pu\`o fare \`e proprio quella che si ottiene \emph{in base al gruppo fondamentale, che serve a definire l'uguaglianza."}}
\end{quote}

The specific purpose of this first paper was to classify linear operators  invariant under coordinate transformations. Though he  lacked a description in terms of group representations and despite not clarifying on which function spaces such operators should act, he was moving towards the same goals as one can find in the contemporary works of Segal (among many others) that we will review in one of the the next sections.

 In this Erlangen--kind of perspective the passage from classical mechanics to Relativity theory can be interpreted as a change in the underlying stability group of the theory. As already explained in \cite{Min}, while classical mechanics is invariant under the action of the $10$--dimensional Galilei group, special relativity is invariant under the group of Poincar\'e transformations\footnote{also referred to as the group of inhomogeneous Lorentz transformations}. The former group can be recovered from the latter as a limit letting a certain parameter $c$, naturally interpreted in physics as the speed of light, go to $+\infty$.

 It is quite natural to ask, at this point, whether this procedure can be iterated, i.e. if the Poincar\'e group can be replaced by other group of symmetries of which it is the limit,  if this extension is unique, if it is further deformable and, finally, when this procedure will stop.
In fact already in \cite[p. 288]{Fa52},   it was remarked that:
\begin{quote}
It may happen that the 10-parameters Lorenz group, which will be denoted by $\mathrm{ Lz}_{3+1}^{10}$, in the evolution of our scientific knowledge, will face an analogous fate as the 10--parameters Galilei group $\mathrm{ Gl}_{3+1}^{10}$ of classical physics and will end up being substituted by a more general group, like the group of isometries $\mathrm{ Ds}_{3+1}^{10}$ of De Sitter space-time. This group, in fact has Lorentz group as limit case, when curvature of space-time goes to $0$, as well as Lorentz group gives as classical limit Galilei group when the speed limit is sent to $\infty$".
\footnote{"Pu\`o accadere che il gruppo di Lorentz, che indicheremo con $\mathrm{Lz}_{3+1}^{10}$ (a 10 parametri) nella progressiva evoluzione delle nostre conoscenze scientifiche, finisca col subire una sorte analoga a quella del gruppo di Galileo $\mathrm{Gl}_{3+1}^{10}$ (pure a 10 parametri) della fisica classica, e cio\`e finisca con l'essere sostituito da un altro gruppo pi\`u generale, per esempio dal gruppo $\mathrm{Ds}_{3+1}^{10}$ dei movimenti in s\`e del cronotopo di De Sitter, che da proprio il gruppo di Lorentz come caso limite, quando il raggio di curvatura del cronotopo tende all'$\infty$, cos\`\i{} come il gruppo di Lorentz da come caso limite il gruppo di Galileo, quando la velocit\`a della luce si fa tendere all'$\infty$."}
\end{quote}

\noindent For this limit procedure he postponed details to a future publication.

This question is, in fact, described by Freeman Dyson as part of one of the \lq\lq missed opportunities"  listed in  his famous Gibbs' Lecture (\cite{Dys}). In his opinion the  reason for such missed opportunities was that physicists and mathematicians were ``neglecting to talk to each other'', and in the specific case Dyson stresses the missed opportunity to foresee an argument in favor of cosmological expansion. To be honest, it would be rather anachronistic to expect an approach like the one of Fantappi\`e at the beginning of the 20th century, as Dyson advocates. De Sitter's work is basically contemporary with the works of Elie Cartan on real semisimple Lie algebras classification (\cite{Car}) and the idea that such Lie algebras could somehow be generic in the set of all Lie algebra laws was yet to come.  The missed opportunity we would rather concentrate on is on how Fantappi\`e lost the chance to anticipate the theory of Lie group contractions and deformations, an idea for which, as we will see in the next paragraphs, times were ripe at the beginning of the 1950's.

 In the 1954 paper (\cite{Fa54}) Fantappi\`e  proved that the inhomogeneous Lorentz group is, truly, a limit of the De Sitter group, which we would nowadays simply call the pseudo-orthogonal group $\mathrm{O}(4,1)$. This proof is the mathematical core of his communication. He then proceeds to outline the basic ideas of  a new relativity theory, which he calls \emph{final relativity}, and to which he will devote some of his residual energies in one of  his last published works  (\cite[p. 158]{Fa55},  emphasis as in the original):

\begin{quote}
And in this Note we will show, indeed, that the Lorentz group is just the \emph{limit case} of \lq\lq another group", continuously depending on a parameter $R$, when $R=\infty$. Since we will show that this other group \emph{cannot be the limit} of a different group, we will call this new group the \emph{final group}, and we will denote it with ${\mathrm Fn}^{10}_{3+1}$.
\footnote{E in questa Nota mostreremo, per l'appunto, che il gruppo di Lorentz \`e proprio il \emph{caso limite} di un \lq\lq altro" gruppo, dipendente con continuit\`a da un parametro $R$, per $R=\infty$, e poich\'e dimostreremo che questo nuovo gruppo \emph{non pu\`o pi\`u essere limite} di un altro gruppo diverso, chiameremo questo nuovo gruppo il \emph{gruppo finale}, e lo indicheremo con ${\mathrm Fn}^{10}_{3+1}$.}
\end{quote} 

It has to be said that, from a physical viewpoint his efforts did not generate results up to his expectations. Despite the fact that his last student, Arcidiacono, devoted most of his scientific career to the development and promotion of Fantappi\`e's final relativity (see \cite{Ar}, and the bibliography therein), his construction was at best ignored and in some case openly criticized. A reflection of the rebuttal by the scientific community can be found in the increasingly sceptical reviews Arcidiacono's papers received on Mathematical Reviews. 

It is not however our aim, here, to judge whether the physical theory of final relativity received the attention it deserved. 
Our point here is to clarify to what extent and why also Fantappi\`e (and his scholars) missed a substantial mathematical opportunity, namely that of pioneering the theory of Lie group and Lie algebra deformations which was being developed at the same time in the USA.

\section{Fantappi\`e's rigidity proof}

In this section we will describe in more detail the content of \cite{Fa54}. In particular we are interested on how Fantappi\`e proved  that the Lorentz group is the limit of a semisimple real Lie group which is not further deformable. This last result, in particular, is what we now term a rigidity result and would be  proved with rather different techniques and considered to be standard, while, at the time, could have led to very interesting new speculations.

We were referring to \lq\lq Lie group deformations" but to say the truth Fantappi\`e's first deformation result concerns, in a way, more the action on a specific homogeneous space than the group itself. Starting from the Lorentz group $\mathrm{ Lz}_{3+1}^{10}$ together with its action as group of isometries of the flat Minkowski space $\mathbb R^4_{3+1}$, he sought for a group $G$ converging, in some limit, to $\lz$ and acting on a manifold $M$. Such manifold was required to fulfill the following:
\begin{enumerate}
\item $\dim M=4$;
\item $G=\mathrm{ Iso}(M)$ has dimension $10$;
\item $M$ is a semi Riemannian manifold of signature $3+1$.
\end{enumerate}

It has to be mentioned that in \cite{Fa54} there is no explicit explanation of the meaning of one group being the limit of another one; it is just 
 intended that there will be a group $G=G(R)$ depending on some parameter $R$ such that when $R=\infty$ then 
$G(\infty)=\lz$. We will further comment on this point later on.

The reasons for asking for the deformed \lq\lq cronotope" $M$ to verify the above requests is explained, in the above mentioned paper at page 159, as follows:
\begin{quote}
Now, the surest informations, valid in the most general case, are clearly those that, referring to \emph{integer numbers} (which cannot vary continuously) should coincide with the known ones, of the limit case
\footnote{Ora, le informazioni pi\`u sicure, che si possono avere nel caso pi\`u generale, sono evidentemente quelle che, riferendosi a \emph{caratteristiche espresse da numeri interi} (le quali non possono variare con continuit\`a) debbono coincidere con quelle note, del caso limite.}.
\end{quote}
This statement, given with no clarification on the involved topologies, may look at first somewhat vague. Though it makes perfect sense, since, no matter which topology is put on some topological space of cronotopes, whatever it may be, functions from such space to $\mathbb Z$ takes constant values on continuous paths. We just remark that no effort is made throughout this paper and subsequent works, on clarifying what a \emph{ space of cronotopes} could be.

The necessity of looking together at $G$ and $M$ (but the same result would hold true without this request) is imposed by the fact that the author chose to rely on a much older result commonly credited to Bianchi \footnote{It is maybe worthwile to recall here that, as mentioned, Bianchi was Fantappi\`e's thesis advisor and most probably his main source of informations for anything concerning differential geometry.}, namely that if $\dim M=4$ and $\dim \mathrm{ Iso}(M)=10$ then $M$ has constant curvature $C$.
Having said so Fantappi\`e utilized the classification of such spaces up to isometries to prove that the requirement on the signature of the semi Riemannian structure enforces $M$ to be a sphere of radius $R$ (where, obviously, $R=1/C^2$) and that in the limit $R\to \infty$ the metric on this spheric cronotope reduces to the usual one on the flat Minkowski space $\mathbb R^4_{3+1}$. At this point he only needed to compute $\mathrm{Iso}(M)$, which he calls the final group $\fn$, and to show that $\fn\simeq \mathrm{O}(4,1)$ as topological (in fact Lie) group\footnote{The choice of giving a special notation like $\fn$ for what is nothing but a pseudo orthogonal group is not completely clear. It is, on one hand, obvious that $\mathrm{Fn}$ stands for final. On the other hand since $\mathrm{Lz}$ stands for Lorentz and $\mathrm{Ga}$ for Galilei, it is difficult not to guess that $\mathrm{Fn}$ could stand also for Fantappi\`e and to ask oneself whether such ambiguity is intended. }.

The second part of the note is centered on the proof of the fact that $\fn$ cannot be further deformed, and is opened by some considerations about $\fn$ being simple. At this point (page 164) the argument goes as follows:
\begin{quote}We can remark that the existence of an invariant subgroup translates into vanishing of some structural constants and, therefore, simplicity translates into the fact that some such constants should be non zero\footnote{An invariant subgroup $H$ in $G$ identifies, via integration, a Lie ideal $\mathfrak h\subseteq\mathfrak g$. By choosing a suitable basis of $\mathfrak g$ one can indeed transform the ideal condition $[x,h]\in\mathfrak h, \forall h\in\mathfrak h, \forall x\in \mathfrak g$ into vanishing of a subset of structural constants with respect to this basis. }. If a simple group is limit of another group, with the same number of parameters, i.e can be obtained by another depending continuously from a variable $\alpha$ for $\alpha\to\alpha_0$, then also its structural constants will depend continuously on $\alpha$, and therefore all those that are non zero for $\alpha=\alpha_0$ will remain $\ne 0$ for $\alpha$ in a neighbourhood of $\alpha_0$. In particular, if the limit group is simple also the variable group should be simple fo $\alpha$ close to $\alpha_0$; in other words a simple group can only be limit of  simple groups.
\footnote{Si pu\`o osservare infatti che la presenza di un sottogruppo invariante si traduce nell'annullarsi di varie costanti di struttura e quindi la \emph{semplicit\`a} del gruppo si traduce nel fatto che almeno alcune di queste costanti debbono essere \emph{diverse da zero}. Se quindi un gruppo semplice \`e \emph{limite} di un altro gruppo, con lo stesso numero di parametri, si ottiene cio\'e da un altro, dipendente con continuit\`a da una variabile $\alpha$ per $\alpha\to \alpha_0$, anche le costanti di struttura di questo dipenderanno con continuit\'a da $\alpha$, e quindi tutte quelle che sono diverse da $0$ per $\alpha=\alpha_0$, resteranno pure $\ne 0$ per $\alpha$ abbastanza prossimo ad $\alpha_0$. In particolare, dunque, se il gruppo limite \`e \emph{semplice}, anche il gruppo variabile deve essere \emph{semplice} per $\alpha$ vicino ad $\alpha_0$, cio\'e \emph{un gruppo semplice non pu\`o essere limite che di un gruppo pure semplice}. }
\end{quote}
This argument is considered to hold for complex simple groups and implies rigidity as a consequence  of complex simple Lie groups classification: there exist no other complex simple Lie groups in dimension $10$. It may be considered of some interest the fact that Fantappi\`e does not refer to the original Cartan's paper but on his later work describing compact forms of complex simple Lie groups (\cite{Ca}). He, then,  looks at the same argument from a real point of view, recalling that $O(5)$, $O(4,1)$ and $O(3,2)$  are the only real simple Lie groups of the right dimension, and remarking that one of these cannot be deformed into another (but about this apparently innocent statement see the quote that ends this section). 

Let us now comment on such proofs.

As for the first part, it has to be noted that the fact that $\lz$ is deformed by $\fn$ is proven by Fantappi\`e's by differential geometric techniques (in a way his deformation takes place in the category of isometry groups of semi--Riemannian manifolds), though it is possible to prove it by purely algebraic techniques. And, in fact, as we shall see in the next paragraph, the existence of such deformation was considered to be quite obvious by contemporaries. Unicity, which is less trivial, was proved by algebraic means by Wigner's student W.T. Sharp in his Ph.D. thesis at Princeton (1960), as credited by Monique L\'evy--Nahas in \cite{LN}. It is difficult to say whether L\'evy--Nahas knew Fantappi\`e's paper (published in Italian on a not so widely known journal) or not. 

As for the proof of rigidity of $\fn$ (which, we recall, is nothing but the real semisimple Lie group $O(4,1)$) his argument is at the same time convincing and specious. Truly he centers the main point: semisemplicity is an open condition and, thus, remains unchanged under small perturbations. However his argument does not meet modern expectations for rigour. The fact that he is talking about structural constants seems to suggest that he is considering Lie algebras rather than Lie groups, a  distinction he never made explicit thus somewhat obscuring his reasoning. In this respect, his work looks closer to papers on Lie theory at the wake of the XXth century, regardless as they are of worries about global considerations. This despite the fact that in the meanwhile, mainly on the basis of H. Weyl's work, topology had entered Lie theory in a substantial way, as clearly explained in \cite{Haw}. While acknowledging that he did make such a distinction in mind, although not clearly expressed, it has to be noticed that the vanishing of some structural constants is strongly dependent, at the Lie algebra level, on the choice of a preferred basis and since he makes no comments on dependence on the choice of this basis, nor on the choice of a canonical one, his argument cannot be regarded as completely conclusive. It has to be remarked that, if correct, this argument would work for a large class of real semisimple Lie groups. Fantappi\`e, however, does not seem to be interested in extending the range of his results beyond this specific case. 

In modern texts, rigidity of real semisimple Lie algebras is proven by vanishing theorems for Lie algebra cohomology. Even if the Whitehead's lemmas were already known at the time (the original Chevalley-Eilenberg paper dates back to 1948) their link with rigidity was not still understood, and would not be for some years  \cite{NR}. Still it is possible that a rigorous proof following Fantappi\`e initial suggestion, i.e. avoiding cohomological arguments, could be given already at the time. It is interesting to note the following remark by Goze and L\"utz remark about the chance of proving rigidity of semisimple Lie algebras by direct algebraic computations in \cite[p. 145]{GoLu}:
\begin{quote}
If you ask some specialists on Lie algebras about this question, he first will tell you that rigidity of semisimple Lie algebras is an evidence, since they are classified: for every $n$ there is only a finite number of non isomorphic semisimple structures on $\mathbb C^n$; hence it is clear, due to continuity, that the only possible limit points are in the class itself. Then you wonder about continuity versus finiteness; naively you think that you could approach $\mu_0$\footnote{$\mu_0$ is the semisimple Lie algebra structure whose rigidity we'd like to show.} jumping from one class to another. You claim for some details... At this point you get the following precisions:
\begin{itemize}
\item[(i)] everybody knows that any Lie algebra with $H^2(\mathfrak g, \mathfrak g)=$ is rigid;
\item[(ii)] it is well-known that semi-simple Lie algebras satisfy $H^2(\mathfrak g,\mathfrak g)=0$, hence are rigid.
\end{itemize}
This continuity versus cohomology is not an evidence for non algebraic-minded people, of course. If you throw an eye in  the literature you discover that (i) is an important result of Nijenhuis and Richardson and that (ii) is a complicated consequence of the existence of Weyl basis using Whitehead Lemmas or alternatively some computations on spectral sequences...On the wole algebraic evidence, isn't it?
\end{quote}

\section{Comparison with  Segal and Inon\"u--Wigner}

We will now consider two related and contemporary works  using approximately the same arguments that are usually credited as the starting point for Lie algebra's deformation (and contraction) theory. It has to be noticed that they both actually predate the work of Fantappi\`e by a few years; assigning priority is not a major concern. It is to be assumed that during and  immediately after the war scientific communication was not at its highs. Switching research topic under such circumstances may quite well result in a limited knowledge of the relevant literature. We have no reason to believe that Fantappi\`e knew such papers, we rather have an indication of the opposite; at Ferrara's library a collection of reprints he owned is preserved and there is no mention of them (\cite{GaPe}). Rather than discussing priorities we will focus on differences between these works and the one of by Fantappi\`e.

\subsection{The work of I.E. Segal on contractions}

 The first paper we consider is by Irving Ezra Segal \cite{Se}, at that time a very young professor in mathematical physics at the University of Chicago. The work was described, in an obituary (\cite{Av}) \emph{a wide ranging article} covering a very large array of algebraic and analytical approaches to quantum field theory.

The focus of this work  is to approach physical systems from an operator algebra point of view, which keeps track of symmetries. Observables are identified inside algebras of operators on Hilbert space (or irreducible representations of $C^*$--algebras). The group of symmetries appears in determining such algebras; observables are generators of unitary representations in bounded or unbounded operators. The largest part of the paper consists of a very careful analysis of a number of technical problems which are encountered when representing locally compact groups (or rather their group $C^*$--algebra) on Hilbert spaces, and when confronting continous and smooth vectors under such representations\footnote{It is no surprise the the theory of smooth vectors in unitary irreps will be definitely set up by Eward Nelson, one of the first PhD's students of Segal.}. On page 255 Segal clearly states that he wants to clarify in what sense  ``one physical theory is a limiting case of the other'', a situation that, in his words, from the point of view of operator algebras ``is quite difficult to define precisely''. This difficulty suggests to him to restrict considerations to Lie algebra deformations. In a very modern approach (which was rediscovered 20 years later) he considers what we nowadays call the algebraic variety of Lie algebra laws, and declares deformations to be continuous paths on such variety. He furthermore proves that any compact Lie algebra is rigid, i.e. it does not admit a nontrivial deformation; in his words \lq\lq the following proof of this fact makes strong use of a suggestion of Harish--Chandra". As for the more general case that includes the non compact real semisimple Lie algebra $\mathfrak{so}(4,1)$ analyzed by Fantappi\`e he says (\cite[p. 257]{Se}):
\begin{quote}
It seems plausible that the result just obtained should be valid without the assumption of compactness, and it can be seen from Cartan's classification of the real simple groups (by the method we use later in connection with the conformal group) that this is true in many particular cases. However, we know of no proof of the general result.
\end{quote}
From this citation it is evident that he was well aware that real semisimple Lie groups' rigidity was an interesting and open problem. More than once in his papers and commentaries in years to follow he expressed his regret for the fact that the starting idea of Lie algebra deformations was credited mostly to the work of In\"on\"u and Wigner, omitting to mention his 1952 paper.
As an aside let us remark that in this work not only does he consider Lie algebra deformations but also more general deformations of a \emph{compound system}, which consists of a locally compact topological space $M$, a continuous action of a locally compact group $G$ on $M$ and a $G$--invariant measure on $M$; a compound system is a topological analogue of the set of data Fantappi\`e is deforming to analyze his \emph{cronotope}, cleared off of the differential data.

In the last part of his paper Segal proceeds to consider a specific situation quite close to the one in \cite{Fa52}. The rigid group of symmetries he ends up with is the $15$--dimensional conformal group $SO(4,2)$. Considering his natural action on the $6$--dimensional space is one way to treat at the same time the De Sitter, anti De Sitter and Poincar\'e group actions on suitable Minkowski space-times, and thus realize Fantappi\`e's deformation as a true limit inside a wider Lie group. He shows that $SO(4,2)$ is not further deformable and sets the foundations of its application to physics, in a program to clarify some aspects of standard special relativity.
It is reasonable to think that since he was able to prove rigidity of $SO(4,2)$ also rigidity of $SO(4,1)$ was clear to him, but no explicit mention of this fact appears in his paper.

Most of Segal's work, in later years, was devoted to the development of a new  relativity theory based on the conformal group (in fact its universal cover), much as Arcidiacono devoted its life to building up in details \emph{final relativity} \cite{Ar}. It is possible that Fantappi\`e's work, and its later developments could be included as a special case in this wide ranging theory. However, as Nelson states (\cite[p. 661]{Av}):
\begin{quote}
Segal's work on the Einstein universe as the arena for cosmology and particle physics is a vast unfinished edifice, constructed with a handful of collaborators.
\end{quote}
A fate it somewhat curiously shared with Fantappi\`e's final relativity.

\subsection{The work of In\"on\"u-Wigner on Lie group contractions}

Let us consider now some comments on \cite{IW52}, a paper which had significant impact on mathematical physics, as evidenced by the 54 citations censed in Mathscinet from 1985 onwards. The introduction (page 510) starts s follows:
\begin{quote}
Classical mechanics is a limiting case of relativistic mechanichs. Hence the group of the former, the Galilei group, must be in some sense a limiting case of the relativistic mechanics' group, the representations of the former must be limiting cases of the latter's representations. There are other examples for similar relations between groups. Thus the inhomogeneous Lorentz group must be, in the same sense, a limiting case of the de Sitter groups\footnote{the plural here refers to the fact that from a purely Lie theory point of view the Poincar\'e group may be reached as a contraction both of $SO(4,1)$ and $SO(3,2)$. A possibility that in \cite{Fa52} is excluded by requirements on the signature of the space-time manifold.}.
\end{quote}
Thus the initial intuition of \cite{Fa52}, for what concerns the existence of the deformation, is stated, immediately, as a fact (as we saw it goes back to Minkowski itself). In fact such an idea was  folklore in US mathematical physics of the time (e.g. ``If this group\footnote{the inhomogeneous Lorentz group} is regarded as the group of motions of a flat three plus one dimensional Riemannian space, it may be regarded as the limit for zero curvature of the group of motions of De Sitter space'' in \cite[p. 113]{Th}). The second remark is that their emphasis is shifted from the beginning from Lie groups to their representations. This is no surprise, since Wigner was one of the pioneers in the use of representation theory techniques  in physics. This makes even more evident how such concept is never explicitly stated in all Fantappi\`e's papers. In fact, a distinction he made in his lectures between the group and its adjoint (see, for example, the posthomous lecture notes of his INDAM course \cite{Fa59}) could be rephrased in terms of left and right action of the same group on a manifold (and on functions on the manifold); a rephrasing which would give a much easier understanding of some  statements.

Going back to In\"on\"u--Wigner the authors proceed to give a definition of Lie groups contraction, which has to be understood as the analogue of the otherwise undefined limit in \cite{Fa54}. It is described in terms of structural constants in the Lie algebra, much closely as by Fantappi\`e, but with no doubts in more exact terms.
Differently from Segal they do not insist on an abstract approach but they give a very clear contraction procedure from an operational point of view. Though such contractions are not as general as Segal's degenerations, as remarked already in the text\footnote{``The above considerations show a certain similarity with those of I.E. Segal. However Segal's considerations are more general than ours as he considers a sequence of Lie groups the structure constants of which converge towards the structure constants of a non isomorphic group. In the above, we have considered only one Lie group but have introduced a sequence of coordinate systems therein and investigated the limiting case of these coordinate systems becoming singular. As a result of our problem being more restricted we could arrive at more specific results.'' \cite[p. 514]{IW52}}, they allow very explicit computations in a number of cases. It has to be said that the first joint work of In\"on\"u and Wigner was published in an Italian journal (\cite{IW52}) exactly in the same years in which Fantappi\`e was developing his theory, and this paper establishes the source of interest in Lie contractions, as explained by the first author in \cite{In}\footnote{This historical note was presented for publication in the proceedings of the mentioned Workshop. However such proceedings were never published. It is possibile to find it on internet at http://ysfine.com/wigner/inonu.pdf.}. However such journal is a journal in physics and it is possible that Fantappi\`e was unaware of  it.
 
On the other hand In\"on\"u and Wigner never discuss the more difficult problem of the uniqueness of deformations; in a way this is reasonable since they take the opposite point of view of contraction. In their approach a reasonable question could have been the following: given a Lie group $G$ list all its possible contractions. Answering such question would requite a complete classification of all its subgroups: a rather ambitious program. Still it is this this work, with its limits, that is nowadays commonly credited as the starting point of Lie algebra deformation theory. It should be no surprise, since Wigner
as Wigner advised many students and coworkers to pursue this avenue thoroughly convinced  of its fruitfulness in mathematical physics.

\section{Conclusions}

The work of Fantappi\`e we discussed may be considered just a partly failed attempt by a mathematician during the end of his career. However, it his clear that his weaknesses, which we emphasized in comparison to the work of Segal and In\"on\"u--Wigner, were not accidental. The \emph{missed opportunity}, in the words of Dyson, was that of foreseeing a fruitful field of study in the direction of Lie algebra deformations and their applications to mathematical physics. His initial intuition was certainly not empty of interesting developments, and his timing was current, as this work was developed at around the same years as the pioneering work of Segal and In\"on\"u--Wigner. 

Even more so if one considers that Fantappi\`e was, at that time, Indam's vice--president and therefore certainly in a position to have a lasting influence on future research directions in his country.

However, several factors conspired to hinder the success of this initial intuition.
First and foremost, of course, Fantappi\`e's own interest was much more centered on pursuing a new \lq\lq final relativity theory" than in understanding the mechanism of symmetries deformation. Two factors, probably, contributed to this. A lack of interest in abstract considerations: this could be influenced by the fact that in the central years of his career abstractness in Italy was negatively connotated as \lq\lq Jewish mathematics" (as explained at various points in \cite{Gu} and, as an example, in the infamous \cite{Ev}). To elaborate further on this point Fantappi\`e was not avoiding generalities: he sought to construct a \lq\lq final relativity" purely on the basis of its symmetries, for example. But he always concentrated on the operational sides of his theories rather than on its foundations, as he did postponing a functional analytic settlement of analytic functionals.

Furthermore he aimed at a complete theory of physics guided by symmetry which, together with its syntropy theory applied to biological systems, would reflect a complete and rational view of the Universe. The idea of deriving such a grand scheme from few  principles followed from his strongly religious views and from his attempts,  more and more intense in the last years of his life, to conciliate faith and science; and effort to which he felt almost compelled (see some personal recallings in \cite{Vt}).

Aside from Fantappi\`e's personal motivation and beliefs it is important to understand the context and work environment; he worked on his \lq\lq final relativity" in relative isolation. No other Italian mathematician would share his interests in Lie theory in the same years, and the global weakness of the Italian school of algebra, especially for what concerned Lie algebras and their representations was still to remain unchanged for many years to come. For what concerns the more physical aspects of his construction  it has to be said that in Italy no one took up the heritage of Levi-Civita, who died in scientific and social isolation, after Ratial Legislation excluded all Jews from teaching in Italian Universities. Thus Fantappi\`e's work, written in Italian, on Italian journals, and devoted to themes which were not much considered by his local contemporaries, remained substantially underdeveloped. It is easy to think that, if health problems had not prevented Fantappi\`e from accepting an intivitation to visit Princeton in 1951 the situation would have been much different. Segal, who was working at  Princeton during the same period, had the chance to interact with Harish--Chandra who was laying the foundations of semisimple Lie groups representation theory. Wigner, who had established his reputation on the use of representation theoretical arguments in particle physics, was as also well positioned for his own students to advance his ideas about group contractions (he met In\"on\"u exactly in 1951). The contrast with the two environments, at that time, couldn't be sharper.

It is sometimes stated that Fantappi\`e's work was neglected due to his beliefs held in religion, philosophy and politics.
 It is certainly possible that some of this more extreme attitudes (e.g. interest on parapsychological phenomena) may have alienated him to an extent. However it is also clear that internal mathematical reasons can be sufficient to explain why his works remained underdeveloped, after his death. More than his personal situation, it was the effect of country's  long period of isolation from the international mathematics community which is why after his death Fantappi\`e's inutitions on Lie group deformations remained one of the \lq\lq missed opportunities" in Italian mathematics research of the time.

\section*{Acknowledgements}

I would like to thank Prof. M.C. Nucci for encouraging me to write these notes and Prof. D. Struppa and Prof. Michael McNeil for generously devoting their time to improve an earlier version of this manuscript.

\end{document}